\documentclass[12pt,english]{article}
\usepackage{theorem}
\usepackage{babel}
\usepackage{a4wide}
\usepackage{graphics}
\usepackage{epsfig}
\usepackage{amssymb}
\usepackage{latexsym}
\usepackage{amsmath}

\newtheorem{thh}{Theorem}[section]

\newtheorem{lem}[thh]{Lemma}
\newtheorem{cor}[thh]{Corollary}
\newtheorem{prop}[thh]{Proposition}

\title{On algebraic automorphisms and \\ their rational invariants}
\author{Philippe Bonnet}
\date{}

\newcommand{\dem}{{\em Proof: }}
\newcommand{\qed}{\begin{flushright} $\blacksquare$\end{flushright}}

\newcommand{\CC}{\mathbb{C}}
\newcommand{\ord}{{\rm{ord}}}
\newcommand{\KB}{\overline{K}}
\newcommand{\CT}{\mathbb{C}^3}
\newcommand{\OX}{{\cal{O}}(X)}
\newcommand{\OXN}{{\cal{O}}(X)^{\nu}}

\begin{document}
\maketitle

\begin{center}
Mathematisches Institut, Universit\"at Basel \\
Rheinsprung 21, 4051 Basel, Switzerland \\ e-mail:
Philippe.bonnet@unibas.ch
\end{center}

\begin{abstract}
Let $X$ be an affine irreducible variety over an algebraically closed field $k$ of characteristic zero.
Given an automorphism $\Phi$, we denote by $k(X)^{\Phi}$ its field of invariants, i.e the set of rational
functions $f$ on $X$ such that $f\circ \Phi=f$. Let $n(\Phi)$ be the transcendence degree of $k(X)^{\Phi}$
over $k$. In this paper, we study the class of automorphisms $\Phi$ of $X$ for which $n(\Phi)=\dim X -1$.
More precisely, we show that under some conditions on $X$, every such automorphism is of the form $\Phi=\varphi_g$, where $\varphi$ is an
algebraic action of a linear algebraic group $G$ of dimension 1 on $X$, and where $g$ belongs to $G$. As an
application, we determine the conjugacy classes of automorphisms of the plane for which $n(\Phi)=1$.
\end{abstract}

\section{Introduction}

Let $k$ be an algebraically closed field of characteristic zero. Let $X$ be an affine irreducible
variety of dimension $n$ over $k$. We denote by $\OX$ its ring of regular functions, and
by $k(X)$ its field of rational functions. Given an algebraic automorphism $\Phi$ of $X$, denote
by $\Phi^*$ the field automorphism induced by $\Phi$ on $k(X)$, i.e. $\Phi^*(f)=f\circ \Phi$ for
any $f\in k(X)$. An element $f$ of $k(X)$ is invariant for $\Phi$ (or simply invariant) if
$\Phi^*(f)=f$. Invariant rational functions form a field denoted $k(X)^{\Phi}$, and we set:
$$
n(\Phi)=trdeg_k \; k(X)^{\Phi}
$$
In this paper, we are going to study the class of automorphisms of $X$ for which $n(\Phi)=n-1$. There
are natural candidates for such automorphisms, such as exponentials of locally nilpotent derivations
(see \cite{M} or \cite{Da}). More generally, one can construct such automorphisms by means of algebraic group
actions as follows. Let $G$ be a linear algebraic group over $k$. An algebraic action of $G$ on $X$ is
a regular map:
$$
\varphi: G\times X \longrightarrow X
$$
of affine varieties, such that $\varphi(g.g',x)=\varphi(g,\varphi(g',x))$ for any $(g,g',x)$ in
$G\times G\times X$. Given an element $g$ of $G$, denote by $\varphi_g$ the
map $x\mapsto \varphi(g,x)$. Then $\varphi_g$ clearly defines an automorphism of $X$. Let $k(X)^G$
be the field of invariants of $G$, i.e. the set of rational functions $f$ on $X$ such that
$f\circ \varphi_g=f$ for any $g\in G$. If $G$ is an algebraic group of dimension $1$, acting
faithfully on $X$, and if $g$ is an element of $G$ of infinite order, then one can prove by
Rosenlicht's Theorem (see \cite{Ro}) that:
$$
n(\varphi_g)=trdeg_k \; k(X)^G= n-1
$$
We are going to see that, under some mild conditions on $X$, there are no other automorphisms
with $n(\Phi)=n-1$ than those constructed above. {\em In what follows, denote by
$\OX^{\nu}$
the normalization of $\OX$, and by $G(X)$ the group of invertible elements of $\OXN$}.

\begin{thh} \label{invariant}
Let $X$ be an affine irreducible variety of dimension $n$ over $k$, such that $char(k)=0$
and $G(X) ^*=k^*$.
Let $\Phi$ be an algebraic automorphism of $X$ such that $n(\Phi)=n-1$. Then there exist an abelian linear
algebraic group $G$ of dimension $1$, and an algebraic action $\varphi$ of $G$ on $X$ such that
$\Phi=\varphi_g$ for some $g\in G$ of infinite order.
\end{thh}
Note that the structure of $G$ is fairly simple. Since every connected linear algebraic group of dimension
1 is either isomorphic to $G_a(k)=(k,+)$ or $G_m(k)=(k^*, \times)$ (see \cite{Hum}, p. 131), there exists a finite abelian
group $H$ such that $G$ is either equal to $H\times G_a(k)$ or $H\times G_m(k)$.
Moreover, the assumption on the group $G(X)$ is essential. Indeed, consider the automorphism
$\Phi$ of $k^* \times k$ given by $\Phi(x,y)=(x,xy)$. Obviously, its field of invariants
is equal to $k(x)$. However, it is easy to check that $\Phi$ cannot have the form
given in the conclusion of Theorem \ref{invariant}.

This theorem is analogous to
a result given by Van den Essen and Peretz (see \cite{V-P}). More precisely, they establish a criterion to decide if
an automorphism $\Phi$ is the exponential of a locally nilpotent derivation, based on the invariants and
on the form of $\Phi$. A similar result has been developed by Daigle (see \cite{Da}).

We apply these results to the group of automorphisms of the plane. First, we obtain a classification of the
automorphisms $\Phi$ of $k^2$ for which $n(\Phi)=1$. Second, we derive a criterion on automorphisms of
$k^2$ to have no nonconstant rational invariants.

\begin{cor} \label{invariant2}
Let $\Phi$ be an algebraic automorphism of $k^2$. If $n(\Phi)=1$, then $\Phi$ is conjugate to one of the
following forms:
\begin{itemize}
\item{$\Phi_1(x,y)=(a^n x,a^m by)$, where $(n,m)\not=(0,0)$, $a,b\in k$, $b$ is a root of unity but $a$ is not,}
\item{$\Phi_2(x,y)=(ax,by+ P(x))$, where $P$ belongs to $k[t]-\{0\}$,  $a,b\in k$ are roots of unity.}
\end{itemize}
\end{cor}

\begin{cor} \label{invariant3}
Let $\Phi$ be an algebraic automorphism of $k^2$. Assume that $\Phi$ has a unique
fixpoint $p$ and that $d\Phi_{p}$ is unipotent. Then $n(\Phi)=0$. 
\end{cor}
We then apply Corollary \ref{invariant3} to an automorphism of $\mathbb{C}^3$
recently discovered by Pierre-Marie Poloni and Lucy Moser (see \cite{M-P}).

We may wonder whether Theorem \ref{invariant} still holds if the ground field $k$ is not algebraically closed or has positive characteristic. The answer is not known for the moment. In fact, two obstructions appear in the proof of Theorem \ref{invariant} when $k$ is arbitrary. First, the group $G_m(k)$ needs to be divisible (see Lemma \ref{Decomposition2}), which is not always the case if $k$ is not algebraically closed. Second, the proof uses the fact that every $G_a(k)$-action on $X$ can be reconstructed from a locally nilpotent derivation on $\OX$ (see subsection \ref{GA}), which is no longer true if $k$ has positive characteristic. This phenomenom is due to the existence of differents forms for the affine line (see \cite{Ru}). Note that, in case Theorem \ref{invariant} holds and $k$ is not algebraically closed, the algebraic group $G$ needs not be isomorphic to $H\times G_a(k)$ or $H\times G_m(k)$, where $H$ is finite. Indeed consider the unit circle $X$ in the plane $\mathbb{R}^2$, given by the equation $x^2 + y^2=1$. Let $\Phi$ be a rotation
in $\mathbb{R}^2$ with center at the origin and angle $\theta\not\in 2\pi \mathbb{Q}$. Then $\Phi$ defines an algebraic automorphism of $X$
with $n(\Phi)=0$, and the subgroup spanned by $\Phi$ is dense in $SO_2(\mathbb{R})$. But $SO_2(\mathbb{R})$ is not isomorphic to either
$G_a(\mathbb{R})$ or $G_m(\mathbb{R})$, even though it is a connected linear algebraic group of dimension 1.

We may also wonder what happens to the automorphisms $\Phi$ of $X$ for which $n(\Phi)=dim X -2$. More precisely, does there exist an action $\varphi$ of a linear algebraic group $G$ on $X$, of dimension 2, such that $\Phi=\varphi_g$ for a given $g\in G$? The answer is no. Indeed consider the automorphism
$\Phi$ of $k^2$ given by $\Phi=f\circ g$, where $f(x,y)=(x+y^2,y)$ and $g(x,y)=(x,y+x^2)$. Let $d(n)$ denote the maximum of the homogeneous degrees of the coordinate functions of the iterate $\Phi^n$. If there existed an action $\varphi$ of a linear algebraic group $G$ such that $\Phi=\varphi_g$, then
the function $d$ would be bounded, which is impossible since $d(n)=4^n$. A similar argument on the length of the iterates also yields the result.
But if we restrict to some specific varieties $X$, for instance $X=k^3$, one may ask the following question: If $n(\Phi)=1$, is $\Phi$ birationally
conjugate to an automorphism that leaves the first coordinate of $k^3$ invariant? The answer is still unknown.

\section{Reduction to an affine curve ${\cal{C}}$}

Let $X$ be an affine irreducible variety of dimension $n$ over $k$. Let $\Phi$ be an algebraic automorphism of $X$ such
that $n(\Phi)=n-1$. In this section, we are going to construct an irreducible affine curve on which $\Phi$ acts naturally.
This will allow us to use some well-known results on automorphisms of curves. We set:
$$
K=\{f \in k(X)|\exists m>0, \; f\circ \Phi^m=f\circ \Phi\circ ...\circ \Phi=f\}
$$
It is straightforward that $K$ is a subfield of $k(X)$ containing both $k$ and $k(X)^{\Phi}$. We begin with some properties
of this field.

\begin{lem} \label{algebrique1}
$K$ has transcendence degree $(n-1)$ over $k$, and is algebraically closed in $k(X)$. In particular, the automorphism $\Phi$
of $X$ has infinite order.
\end{lem}
\dem First we show that $K$ has transcendence degree $(n-1)$ over $k$. Since $K$ contains the field $k(X)^{\Phi}$, whose
transcendence degree is $(n-1)$, we only need to show that the extension $K/k(X)^{\Phi}$ is algebraic, or in other words
that every element of $K$ is algebraic over $k(X)^{\Phi}$. Let $f$ be any element of $K$. By definition, there exists an
integer $m>0$ such that $f\circ \Phi^m=f$. Let $P(t)$ be the polynomial of $k(X)[t]$ defined as:
$$
P(t)=\prod_{i=0} ^{m-1} (t - f\circ \Phi^i)
$$
By construction, the coefficients of this polynomial are all invariant for $\Phi$, and $P(t)$ belongs to $k(X)^{\Phi}[t]$.
Moreover $P(f)=0$, $f$ is algebraic over $k(X)^{\Phi}$ and the first assertion follows.

Second we show that $K$ is algebraically closed in $k(X)$. Let $f$ be an element of $k(X)$ that is algebraic over $K$.
We need to prove that $f$ belongs to $K$. By the first assertion of the lemma, $f$ is algebraic over $k(X)^{\Phi}$. Let
$P(t)=a_0 + a_1t + ...+ a_pt^p$ be a nonzero
minimal polynomial of $f$ over $k(X)^{\Phi}$. Since $P(f)=0$ and all $a_i$ are invariant, we have
$P(f\circ \Phi)=P(f)\circ \Phi=0$. In particular, all elements of the form $f\circ \Phi^i$, with $i\in \mathbb{N}$,
are roots of $P$. Since $P$ has finitely many roots, there exist two distinct integers $m'<m''$ such that
$f\circ \Phi^{m'}=f\circ \Phi^{m''}$. In particular, $f\circ \Phi^{m''-m'}=f$ and $f$ belongs to $K$.

Now if $\Phi$ were an automorphism of finite order, then $K$ would be equal to $k(X)$. But this is impossible since
$K$ and $k(X)$ have different transcendence degrees.
\qed

\begin{lem} \label{algebrique2}
There exists an integer $m>0$ such that $K=k(X)^{\Phi^m}$.
\end{lem}
\dem By definition, $k(X)$ is a field of finite type over $k$. Since $K$ is contained in $k(X)$, $K$ has also finite
type over $k$. Let $f_1,...,f_r$ be some elements of $k(X)$ such that $K=k(f_1,...,f_r)$. Let $m_1,...,m_r$ be some
positive integers such that $f_i\circ\Phi^{m_i}=f_i$, and set $m=m_1...m_r$. By construction, all $f_i$ are invariant
for $\Phi^m$. In particular, $K$ is invariant for $\Phi^m$ and $K\subseteq k(X)^{\Phi^m}$. Since $k(X)^{\Phi^m}\subseteq K$,
the result follows.
\qed 
Let $L$ be the algebraic closure of $k(X)$, and let $A$ be the $K$-subalgebra of $L$ spanned by $\OX$. By construction,
$A$ is an integral $K$-algebra of finite type of dimension 1. Let $m$ be an integer satisfying the conditions of lemma
\ref{algebrique2}. The automorphism $\Psi^*=(\Phi^m)^*$ of $\OX$ stabilizes $A$, hence it defines
a $K$-automorphism of $A$, of infinite order (see lemma \ref{algebrique1}). Let $B$ be the integral closure
of $A$. Then $B$ is also an integral $K$-algebra of finite type, of dimension 1, and the $K$-automorphism
$\Psi^*$ extends uniquely to $B$. If $\overline{K}$ stands for the algebraic closure of $K$, we set:
$$
C=B\otimes_K \overline{K}
$$
By construction, ${\cal{C}}=Spec(C)$ is an affine curve over the algebraically closed field $\overline{K}$. Moreover
the automorphism $\Psi^*$ acts on $C$ via the operation:
$$
\Psi^*: C\longrightarrow C, \quad x\otimes y \longmapsto \Psi^*(x)\otimes y
$$
This makes sense since $\Psi^*$ fixes the field $K$. Therefore $\Psi^*$ induces an algebraic automorphism of the
curve ${\cal{C}}$. Since $K$ is algebraically closed in $k(X)$ by lemma \ref{algebrique1}, $C$ is integral
(see \cite{Z-S}, Chap. VII, \S 11, Theorem 38). But by construction, $B$ and $\overline{K}$ are normal rings. Since $C$ is a domain and $char(K)=0$, $C$ is
also integrally closed by a result of Bourbaki (see \cite{Bou}, p. 29). So $C$ is a normal domain
and ${\cal{C}}$ is a smooth irreducible curve. 

\begin{lem} \label{Courbe}
Let $C$ be the $\KB$-algebra constructed above. Then either $C=\overline{K}[t]$ or $C=\overline{K}[t,1/t]$.
\end{lem}
\dem By lemma \ref{algebrique1}, the automorphism $\Phi$ of $X$ has infinite order. Since the fraction field of $B$ is equal to
$k(X)$, $\Psi^*$ has infinite order on $B$. But $B\otimes 1\subset C$, so $\Psi^*$ has infinite order on $C$. In particular,
$\Psi$ acts like an automorphism of infinite order on ${\cal{C}}$. Since ${\cal{C}}$ is affine, it has genus zero (see \cite{Ro2}).
Since $\KB$ is algebraically closed, the curve ${\cal{C}}$ is rational (see \cite{Che}, p. 23 ). Since ${\cal{C}}$ is smooth, it is isomorphic
to $\mathbb{P}^1(\KB)-E$, where $E$ is a finite set. Moreover, $\Psi$ acts like an automorphism of $\mathbb{P}^1(\KB)$ that stabilizes $\mathbb{P}^1(\KB)-E$. Up to replacing
$\Psi$ by one of its iterates, we may assume that $\Psi$ fixes every point of $E$. But an automorphism of $\mathbb{P}^1(\KB)$ that fixes at least three points
is the identity, which is impossible. Therefore $E$ consists of at most two points, and ${\cal{C}}$ is either isomorphic to $\KB$ or
to $\KB^*$. In particular, either $C=\overline{K}[t]$ or $C=\overline{K}[t,1/t]$.
\qed 

\section{Normal forms for the automorphism $\Psi$}

Let $C$ and $\Psi^*$ be the $\KB$-algebra and the $\KB$-automorphism constructed in the previous section. In this section, we are going to give normal forms for the couple $(C,\Psi^*)$, in case the group
$G(X)$ is trivial, i.e. $G(X)=k^*$. We begin with a few lemmas.

\begin{lem} \label{Diviseur}
Let $X$ be an irreducible affine variety over $k$. Let $\Psi$ be an automorphism of $X$. Let
$\alpha,f$ be some elements of $k(X)^*$ such that $(\Psi^*)^n(f)=\alpha ^n f$ for any $n\in \mathbb{Z}$. Then $\alpha$ belongs to $G(X)$.
\end{lem}
\dem Given an element $h$ of $k(X)^*$ and a prime divisor $D$ on the normalization $X^{\nu}$, we consider $h$ as a rational function on $X^{\nu}$, and denote by $\ord_{D}(h)$ the multiplicity of $h$ along $D$.
This makes sense since the variety $X^{\nu}$ is
normal. Fix any prime divisor $D$ on $X$. Since $(\Psi^*)^n(f)=\alpha ^n f$ for any
$n\in \mathbb{Z}$, we obtain:
$$
\ord_{D}((\Psi^*)^n(f))= n\ord_{D}(\alpha) + \ord_D(f)
$$
Since $\Psi$ is an algebraic automorphism of $X$, it extends uniquely to an algebraic automorphism of
$X^{\nu}$, which is still denoted $\Psi$. Moreover, this extension maps every prime divisor to another prime divisor, does not change the multiplicity
and maps distinct prime divisors into distinct ones. If $div(f)=\sum_i n_i D_i$, where all $D_i$ are prime, then we have:
$$
div((\Psi^*)^n(f))=\sum_i n_i (\Psi^*)^n(D_i)
$$
where all $(\Psi^*)^n(D_i)$ are prime and distinct. So the multiplicity of $(\Psi^*)^n(f)$ along $D$ is equal to zero
if $D$ is none of the $(\Psi^*)^n(D_i)$, and equal to $n_i$ if $D=(\Psi^*)^n(D_i)$. In all cases, if $R=max\{|n_i|\}$, then we find
that $|\ord_{D}((\Psi^*)^n(f))|\leq R$ and $|\ord_D(f)|\leq R$, and this implies for any integer $n$:
$$
|n\ord_{D}(\alpha)|\leq 2R
$$
In particular we find $\ord_{D}(\alpha)=0$. Since this holds for any prime divisor $D$, the support of $div(\alpha)$ in $X^{\nu}$ is empty and $div(\alpha)=0$. Since $X^{\nu}$ is normal, $\alpha$ is an invertible element of $\OXN$, hence it belongs to $G(X)$.
\qed

\begin{lem} \label{Diviseur2}
Let $K$ be a field of characteristic zero and $\KB$ its algebraic closure. Let $C$ be either equal to $\KB[t]$ or to $\KB[t,1/t]$. Let
$\Psi^*$ be a $\KB$-automorphism of $C$ such that $\Psi^*(t)=at$, where $a$ belongs to $\KB$. Let $\sigma_1$ be a $K$-automorphism of $C$, commuting
with $\Psi^*$, such that $\sigma_1(\KB)=\KB$. Then $\sigma_1(a)$ is either equal to $a$ or to $1/a$.
\end{lem}
\dem We distinguish two cases depending on the ring $C$. First assume that $C=\KB[t]$. Since $\sigma_1$ is a $K$-automorphism of $C$ that
maps $\KB$ to itself, we have $\KB[t]=\KB[\sigma_1(t)]$. In particular $\sigma_1(t)= \lambda t + \mu$, where $\lambda,\mu$ belong to $\KB$
and $\lambda\not=0$. Since $\Psi^*$ and $\sigma_1$ commute, we obtain:
$$
\Psi^* \circ \sigma_1(t)=\lambda at +\mu = \sigma_1 \circ \Psi^*(t)=\sigma_1(a)(\lambda t +\mu)
$$
In particular, we have $\sigma_1(a)=a$ and the lemma follows in this case. Second assume that $C=\KB[t,1/t]$. Since $\sigma_1$ is a
$K$-automorphism of $C$, we find:
$$
\sigma_1(t)\sigma_1(1/t)=\sigma_1(t.1/t)=\sigma_1(1)=1
$$
Therefore $\sigma_1(t)$ is an invertible element of $C$, and has the form $\sigma_1(t)=a_1 t^{n_1}$, where $a_1\in \KB^*$ and $n_1$
is an integer. Since $\sigma_1$ is a $K$-automorphism of $C$ that maps $\KB$ to $\KB$, we have $\KB[t,1/t]=\KB[\sigma_1(t),1/\sigma_1(t)]$.
In particular $|n_1|=1$ and either $\sigma_1(t)=a_1t$ or $\sigma_1(t)=a_1/t$. If $\sigma_1(t)=a_1t$, the relation $\Psi^*\circ\sigma_1(t)
=\sigma_1\circ \Psi^*(t)$ yields $\sigma(a)=a$. If $\sigma_1(t)=a_1/t$, then the same relation yields $\sigma(a)=1/a$.
\qed

\begin{lem} \label{Inversible2}
Let $X$ be an irreducible affine variety of dimension $n$ over $k$, such that $G(X)=k^*$. Let $\Phi$ be an automorphism of $X$ such
that $n(\Phi)=(n-1)$. Let $\Psi^*$ be the automorphism of $C$ constructed in the previous section. If
either $C=\KB[t]$ or $C=\KB[t,1/t]$, and if $\Psi^*(t)=at$, then $a$ belongs to $k^*$.
\end{lem}
\dem We are going to prove by contradiction that $a$ belongs to $k^*$. So assume that $a\not\in k^*$.
Let $\sigma$ be any element of $Gal(\KB/K)$, and denote by $\sigma_1$ the $K$-automorphism of
$C$ defined as follows:
$$
\forall (x,y)\in B\times \KB, \; \sigma_1(x\otimes y)=x\otimes\sigma_1(y)
$$
Since $\Psi^* \circ \sigma_1 (x\otimes y)=\Psi^*(x)\otimes \sigma_1(y)=\sigma_1 \circ \Psi^*(x\otimes y)$ for any element $x\otimes y$ of $B\otimes_K \KB$, $\Psi^*$
and $\sigma_1$ commute. Moreover if we identify $\KB$ with $1\otimes \KB$, then $\sigma_1(\KB)=\KB$ by construction. By lemma \ref{Diviseur2},
we obtain:
$$
\forall \sigma \in Gal(\KB/K), \quad \sigma(a)=a \quad \mbox{or} \quad \sigma(a)=\frac{1}{a}
$$
In particular, the element $(a^i+a^{-i})$ is invariant under the action of $Gal(\KB/K)$ for any $i$, and so it belongs to $K$ because $char(K)=0$.
Now let $f$ be an element of $B-K$. Since $f$ belongs to $C$, we can express $f$ as follows:
$$
f=\sum_{i=r} ^s f_i t^i
$$
Choose an $f\in B-K$ such that the difference $(s-r)$ is minimal. We claim that $(s-r)=0$, i.e. $f=f_s t^s$. Indeed, assume that $s>r$.
Since $f$ is an element of $B$, the following expressions:
$$
\begin{array}{rcl}
\Psi^*(f) + (\Psi^*)^{-1}(f) - (a^s+a^{-s})f &=&\sum_{i=r} ^{s-1}  f_i (a^i + a^{-i} -a^s -a^{-s})t^i \\
\Psi^*(f) + (\Psi^*)^{-1}(f) - (a^r+a^{-r})f &=&\sum_{i=r+1} ^{s}  f_i (a^i + a^{-i} -a^r -a^{-r})t^i
\end{array}
$$
also belong to $B$. By minimality of $(s-r)$, these expressions belong to $K$. In other words, $f_i (a^i + a^{-i} -a^s -a^{-s})=0$
(resp. $f_i (a^i + a^{-i} -a^r -a^{-r})=0$) for any $i\not=0,s$ (resp. for any $i\not=0,r$). Since $k$ is algebraically closed and
$a\not\in k^*$ by assumption, $(a^i + a^{-i} -a^s -a^{-s})$ (resp. $(a^i + a^{-i} -a^r -a^{-r})$) is nonzero for any $i\not=s$ (resp.
for any $i\not=r$). Therefore $f_i=0$ for any $i\not=0$, and $f$ belongs to $K$, a contradiction. Therefore $s=r$ and $f=f_st^s$.
Since $f$ belongs to $B$, it also belongs to $k(X)$. Since $\Psi$ is an automorphism of $X$, the element $a^s=\Psi^*(f)/f$ belongs
to $k(X)$. Moreover $(\Psi^*)^n(f)=a^{ns}f$ for any $n\in \mathbb{Z}$. By lemma \ref{Diviseur}, $a^s$ belongs to $G(X)=k^*$. Since $k$ is algebraically closed, $a$ belongs to $k^*$, hence a contradiction, and the result follows.
\qed

\begin{prop} \label{Courbe2}
Let $X$ be an irreducible affine variety of dimension $n$ over $k$, such that $G(X)=k^*$. Let $\Phi$ be an automorphism of $X$ such
that $n(\Phi)=(n-1)$. Let $C$ and $\Psi^*$ be the $\KB$-algebra and the $\KB$-automorphism constructed in the previous section. Then
up to conjugation, one of the following three cases occurs:
\begin{itemize}
\item{$C=\KB[t]$ and $\Psi^*(t)=t+1$,}
\item{$C=\KB[t]$ and $\Psi^*(t)=at$, where $a\in k^*$ is not a root of unity,}
\item{$C=\KB[t,1/t]$ and $\Psi^*(t)=at$, where $a\in k^*$ is not a root of unity.}
\end{itemize}
\end{prop}
\dem By lemma \ref{Courbe}, we know that either $C=\overline{K}[t]$ or $C=\overline{K}[t,1/t]$. We are going to study both cases.\\ \\
{\underline{First case}}: $C=\KB[t]$. \\ \\
The automorphism $\Psi^*$ maps $t$ to $at+b$, where $a\in \KB^*$ and $b\in \KB$. If $a=1$, then $b\not=0$ and up to replacing $t$ with $t/b$,
we may assume that $\Psi^*(t)=t+1$. If $a\not=1$, then up to replacing $t$ with
$t-c$ for a suitable $c$, we may assume that $\Psi^*(t)=at$. But then lemma \ref{Inversible2} implies that $a$ belongs to $k^*$. Since
$\Psi^*$ has infinite order, $a$ cannot be a root of unity. \\ \\
{\underline{Second case}}: $C=\KB[t,1/t]$. \\ \\
Since $\Psi^*(t)\Psi^*(1/t)=\Psi^*(1)=1$, $\Psi^*(t)$ is an invertible element of $C$. So $\Psi^*(t)=at^n$, where $a\in \KB^*$ and $n\not=0$.
Since $\Psi^*$ is an automorphism, $n$ is either equal to $1$ or to $-1$. But if $n$ were equal to $-1$, then a simple computation
shows that $(\Psi^*)^2$ would be the identity, which is impossible. So $\Psi^*(t)=at$, where $a\in \KB^*$. By lemma \ref{Inversible2}, $a$
belongs to $k^*$. As before, $a$ cannot be a root of unity.
\qed

\section{Proof of the main theorem}

In this section, we are going to establish Theorem \ref{invariant}. We will split its proof in two steps
depending on the form of the automorphism $\Psi^*$ given in Proposition \ref{Courbe2}. But before, we begin
with a few lemmas.

\begin{lem} \label{Commute1}
Let $\Phi$ be an automorphism of an affine irreducible variety $X$. Let $G$ be a linear algebraic group and
$\psi$ be an algebraic $G$-action on $X$. Let $h$ be an element of $G$ such that the group $<\!h\!>$
spanned by $h$ is Zariski dense in $G$. If $\Phi$ and $\psi_h$ commute, then $\Phi$ and $\psi_g$
commute for any $g$ in $G$.
\end{lem}
\dem It suffices to check that $\Phi^*$ and $\psi^* _g$ commute for any $g\in G$. For any $k$-algebra
automorphisms $\alpha,\beta$ of $\OX$, denote by $[\alpha,\beta]$ their commutator, i.e.
$[\alpha,\beta]=\alpha \circ\beta \circ \alpha^{-1} \circ \beta^{-1}$. For any $f\in \OX$, set:
$$
\lambda(g,f)(x)=[\Phi^*,\psi^* _g ](f)(x)-f(x)
$$
Since $G$ is a linear algebraic group acting algebraically on the affine variety $X$, $\lambda(g,f)(x)$
is a regular function
on ${G\times X}$. Since $\Phi^*$ and $\psi^* _h$ commute,
the automorphisms $\Phi^*$ and $\psi^* _{h^n}$ commute for any integer $n$. So the regular function
$\lambda(g,f)(x)$ vanishes on $<\!h\!>\times X$. Since $<\!h\!>$ is dense in $G$ by
assumption, $<\!h\!>\times X$ is dense in $G\times X$ and $\lambda(g,f)(x)$ vanishes identically
on $G\times X$ . In particular, $[\Phi^*,\psi^* _g](f)=f$ for any $g\in G$.
Since this holds for any element $f$ of $\OX$, the bracket $[\Phi^*,\psi^* _g]$
coincides with the identity on $\OX$ for any $g\in G$, and the result follows.
\qed

\begin{lem} \label{Decomposition2}
Let $\Phi$ be an automorphism of an affine irreducible variety $X$. Let $G$ be a linear
algebraic group
and $\psi$ be an algebraic $G$-action on $X$. Let $h$ be an element of $G$ such that the group $<\!h\!>$
spanned by $h$ is Zariski dense in $G$. Assume there exists a nonzero integer $r$ such that $\Phi^r = \psi_h$,
and that $G$ is divisible.
Then there exists an algebraic action $\varphi$ of $G'= \mathbb{Z}/r\mathbb{Z}
\times G$ such that $\Phi=\varphi_{g'}$ for some $g'$ in $G'$.
\end{lem}
\dem Fix an element $b$ in $G$ such that
$b^r=h$, and set $\Delta=\Phi \circ \psi_{b^{-1}}$. This is possible since $G$ is divisible.
By construction, $\Delta$ is an automorphism of $X$.
Since $\Phi^r= \psi_h$, $\Phi$ and $\psi_h$ commute. By lemma \ref{Commute1}, $\Phi$ and $\psi_g$ commute
for any $g\in G$. In particular, we have:
$$
\Delta^r= (\Phi^r)\circ  \psi_{b^{-r}}=(\Phi^r) \circ \psi_{h^{-1}}=\mbox{Id}
$$ 
So $\Delta$ is finite, $\Phi = \Delta \circ \psi_{b}$ and $\Delta$ commutes with $\psi_g$
for any $g \in G$. The group $G'$ then acts on $X$ via the map $\varphi$ defined by:
$$
\varphi_{(i,g)}(x)=\Delta^i \circ \psi_g(x)
$$
Moreover we have $\Phi=\varphi_{g'}$ for $g'=(1,b)$.
\qed
The proof of Theorem \ref{invariant} will then go as follows. In the following subsections, we are going
to exhibit an algebraic action $\psi$ of $G_a(k)$ (resp. $G_m(k)$) on $X$, such that
$\Psi=\Phi^m=\psi_h$ for some $h$. In both cases, the group $G$ we will consider will be linear
algebraic of dimension 1, and divisible. Moreover the element $h$ will span a Zariski dense
set because $h\not=0$ (resp. $h$ is not a root of unity). With these conditions, Theorem
\ref{invariant} will become a direct application of Lemma \ref{Decomposition2}.

\subsection{The case $\Psi^*(t)=t+1$} \label{GA}

Assume that $C=\KB[t]$ and $\Psi^*(t)=t+1$. We are going
to construct a nontrivial algebraic $G_a(k)$-action $\psi$ on $X$ such that $\Psi=\psi_1$.
Since $\OX\subset C$,
every element $f$ of $\OX$ can be written as $f=P(t)$, where $P$ belongs
to $\KB[t]$. We set $r=\deg_t P(t)$. Since $\Psi^*$ stabilizes $\OX$, the expression:
$$
(\Psi^i)^*(f)=P(t+i)=\sum_{j=0} ^r P^{(j)}(t)\frac{i^j}{j!}
$$
belongs to $\OX$ for any integer $i$. Since the matrix $M=(i^j/j!)_{0\leq i,j\leq r}$ is invertible in
${\cal{M}}_{r+1}(\mathbb{Q})$, the polynomial $P^{(j)}(t)$ belongs to $\OX$ for any $j\leq r$. So the
$\KB$-derivation $D=\partial/\partial t$ on $C$ stabilizes the $k$-algebra $\OX$. Since $D^{r+1}(f)=0$,
the operator $D$, considered as a $k$-derivation on $\OX$, is locally nilpotent (see \cite{Van}).
Therefore the exponential map:
$$
\exp uD: \OX \longrightarrow \OX[u], \quad f \longmapsto \sum_{j\geq 0} D^j(f)\frac{u^j}{j!}
$$
is a well-defined $k$-algebra morphism. But $\exp uD$ defines also a $K$-algebra morphism from $C$
to $C[u]$. Since $\exp uD(t)=t+u$, $\exp D$ coincides with $\Psi^*$ on $C$. Since $C$ contains the
ring $\OX$, we have $\exp D=\Psi^*$ on $\OX$. So the exponential map induces an algebraic
$G_a(k)$-action $\psi$ on $X$ such that $\Psi= \psi_1$ (see \cite{Van}).

\subsection{The case $\Psi^*(t)=at$}

Assume that $\Psi^*(t)=at$ and $a$ is not
a root of unity. We are going to construct a nontrivial algebraic $G_m(k)$-action $\psi$ on $X$
such that $\Psi=\psi_a$. First note that either $C=\KB[t]$ or $C=\KB[t,1/t]$. Let $f$ be any element
of $\OX$. Since $\OX\subset C$, we can write $f$ as:
$$
f=P(t)=\sum_{i=r} ^s f_i t^i
$$
where the $f_it^i$ belong a priori to $C$. Since $\Psi^*$ stabilizes $\OX$, the expression:
$$
(\Psi^j)^*(f)=P(a^jt)=\sum_{i=r} ^s a^{ji} f_i t^i
$$
belongs to $\OX$ for any integer $j$. Since $a$ belongs to $k^*$ and is not a root of unity, the Vandermonde
matrix $M=(a^{ij})_{0\leq i,j\leq s-r}$ is invertible in ${\cal{M}}_{s-r+1}(k)$. So the elements
$f_it^i$ all belong to $\OX$ for any integer $i$. Consider the map:
$$
\psi^*: \OX \longrightarrow \OX[v,1/v], \quad f\longmapsto \sum_{i=r} ^s f_i t^i v^i
$$
Then $\psi^*$ is a well-defined $k$-algebra morphism, which induces a regular map $\psi$ from $k^*\times X$ to $X$.
Moreover we have $\psi_v\circ \psi_{v'}=\psi_{vv'}$ on $X$ for any $v,v'\in k^*$. So $\psi$ defines an algebraic
$G_m(k)$-action on $X$ such that $\Psi=\psi_a$.

\section{Proof of Corollary \ref{invariant2}}

Let $\Phi$ be an automorphism of the affine plane $k^2$, such that $n(\Phi)=1$. By Theorem \ref{invariant},
there exists an algebraic action $\varphi$ of an abelian linear algebraic group $G$ of dimension $1$
such that $\Phi=\varphi_g$. We will distinguish the cases
$G=\mathbb{Z}/r\mathbb{Z} \times G_m(k)$ and $G=\mathbb{Z}/r\mathbb{Z} \times G_a(k)$. \\ \\
{\underline{First case}}: $G=\mathbb{Z}/r\mathbb{Z} \times G_m(k)$. \\ \\
Then $G$ is linearly reductive and $\varphi$ is conjugate to a representation in $Gl_2(k)$ (see \cite{Ka}
or \cite{Kr}). Since $G$ consists solely of semisimple elements,
$\varphi$ is even diagonalizable. In particular, there exists a system $(x,y)$ of polynomial
coordinates, some integers $n,m$ and some $r$-roots of unity $a,b$ such that:
$$
\varphi_{(i,u)}(x,y)=(a^i u^{n} x, b^i u^m y)
$$
Note that, since the action is faithful, the couple $(n,m)$ is distinct from $(0,0)$. Since
$k$ is algebraically closed, we can even reduce $\Phi=\varphi_g$ to the first form given
in Corollary \ref{invariant2}. \\ \\
{\underline{Second case}}: $G=\mathbb{Z}/r\mathbb{Z} \times G_a(k)$. \\ \\
Let $\psi$ and $\Delta$ be respectively the $G_a(k)$-action and finite automorphism constructed
in Lemma \ref{Decomposition2}. By Rentschler's theorem (see \cite{Re}), there exists a system
$(x,y)$ of polynomial coordinates and an element $P$ of $k[t]$ such that:
$$
\psi_u(x,y)=(x,y +uP(x))
$$
For any $f\in k[x,y]$, set $\deg_{\psi}(f)=\deg_u \exp uD(f)$. It is well-known that this defines
a degree function on $k[x,y]$ (see \cite{Da}). Since $\psi$ and $\Delta$ commute, $\Delta^*$ preserves the space
$E_n$ of polynomials of degree $\leq n$ with respect to $\deg_{\psi}$. In particular,
$\Delta^*$ preserves $E_0=k[x]$. So $\Delta^*$ induces a finite automorphism of $k[x]$,
hence $\Delta^*(x)=ax+b$, where $a$ is a root of unity. Since $\Delta$ is finite, either
$a\not=1$ or $a=1$ and $b=0$. In any case, up to replacing $x$ by $x-\mu$ for a suitable
constant $\mu$, we may assume that $\Delta^*(x)=ax$. Moreover $\Delta^*$ preserves
the space $E_1=k[x]\{1,y\}$. With the same arguments as before, we obtain that
$\Delta^*(y)=cy+d(x)$, where $c$ is a root of unity and $d(x)$ belongs to $k[x]$.
Composing $\Delta$ with $\psi_{1/m}$ then yields the second form given
in Corollary \ref{invariant2}.

\section{Proof of Corollary \ref{invariant3}}

Let $\Phi$ be an algebraic automorphism of $k^2$. We assume that $\Phi$ has a unique
fixpoint $p$ and that $d\Phi_{p}$ is unipotent. We are going to prove that $n(\Phi)=0$.

First we check that $n(\Phi)$ cannot be equal to $2$. Assume that $n(\Phi)=2$. Then $k(x,y)^{\Phi}$ has
transcendence degree 2, and the extension $k(x,y)/k(x,y)^{\Phi}$ is algebraic, hence finite. Moreover $\Phi^*$
acts like an element of the Galois group of this extension. In particular, $\Phi^*$ is finite. By a result of Kambayashi
(see \cite{Ka}), $\Phi$ can be written as $h\circ A\circ h^{-1}$, where $A$ is an element
of $Gl_2(k)$ of finite order and $h$ belongs to $Aut(k^2)$. Since $\Phi$ has a unique
fixpoint $p$, we have $h(0,0)=p$. In particular, $d\Phi_{p}$ is conjugate to $A$ in $Gl_2(k)$.
Since $d\Phi_{p}$ is unipotent and $A$ is finite, $A$ is the identity. Therefore $\Phi$
is also the identity, which contradicts the fact that it has a unique fixpoint. 

Second we check that $n(\Phi)$ cannot be equal to $1$. Assume that $n(\Phi)=1$. By the previous
corollary, up to conjugacy, we may assume that $\Phi$ has one of the following forms:
\begin{itemize}
\item{$\Phi_1(x,y)=(a^n x,a^m by)$, where $(n,m)\not=(0,0)$, $b$ is a root of unity but $a$ is not,}
\item{$\Phi_2(x,y)=(ax,by+ P(x))$, where $P$ belongs to $k[t]-\{0\}$ and $a,b$ are roots of unity.}
\end{itemize}
Assume that $\Phi$ is an automorphism of type $\Phi_1$. Then $d\Phi_{p}$ is a diagonal matrix
of $Gl_2(k)$, distinct from the identity. But this is impossible since $d\Phi_{p}$ is unipotent.
So assume that $\Phi$ is an automorphism of type $\Phi_2$. Then $d\Phi_{p}$ is a linear
map of the form $(u,v)\mapsto (au,bv +du)$, with $d\in k$. Since $d\Phi_{p}$ is unipotent, we
have $a=b=1$. So $(\alpha,\beta)$ is a fixpoint if and only if $P(\alpha)=0$. In particular, the set
of fixpoints is either empty or a finite union of parallel lines. But this is impossible
since there is only one fixpoint by assumption. Therefore $n(\Phi)=0$.

\section{An application of Corollary \ref{invariant3}}

In this section, we are going to see how Corollary \ref{invariant3} can be applied to the determination
of invariants for automorphisms of $\mathbb{C}^3$. Set $Q(x,y,z)=x^2 y - z^2 -xz^3$ and consider the
following automorphism (see \cite{M-P}):
$$
\Phi: \CT\longrightarrow \CT, \quad 
\begin{pmatrix}
x \\
y \\
z
\end{pmatrix}
\longmapsto
\begin{pmatrix}\
x \\
y(1-xz) +\frac{Q^2}{4} + z^4\\
z - \frac{Q}{2} x
\end{pmatrix}
$$
We are going to show that:
$$
\mathbb{C}(x,y,z)^{\Phi}=\mathbb{C}(x) \quad \mbox{and} \quad \mathbb{C}[x,y,z]^{\Phi}=\mathbb{C}[x]
$$
Let $k$ be the algebraic closure of $\mathbb{C}(x)$. Since $\Phi^*(x)=x$, the morphism $\Phi^*$ induces an
automorphism of $k[y,z]$, which we denote by $\Psi^*$. The automorphism $\Psi$ has clearly $(0,0)$ as a fixpoint,
and its differential at this point is unipotent, distinct from the identity. Indeed, it is given by the matrix:
$$
d\Psi_{(0,0)}=
\begin{pmatrix}
1 & -x^3/2 \\
0 & 1
\end{pmatrix}
$$
Moreover, the set of fixpoints of $\Psi$ is reduced to the origin. Indeed, if $(\alpha,\beta)$ is a point of $k^2$
fixed by $\Psi$, then $xQ=0$ and $4\beta^4 -4x\alpha\beta +Q^2=0$. Since $x$ belongs to $k^*$, we have:
$$
Q=x^2 \alpha -\beta^2 -x\beta^3=0 \quad \mbox{and} \quad \beta^4 -x\alpha\beta=0
$$
If $\beta=0$, then $\alpha=0$ and we find the origin. If $\beta\not=0$, then dividing by $\beta$ and multiplying by $-x$ yields the relation:
$$
x^2\alpha - x\beta^3=0
$$
This implies $\beta^2=0$ and $\beta=0$, hence a contradiction. By Corollary \ref{invariant3}, the field of invariants
of $\Psi$ has transcendence degree zero. So the field of invariants of $\Phi$ has transcendence degree $\leq 1$
over $\mathbb{C}$. Since this field contains $\mathbb{C}(x)$ and that $\mathbb{C}(x)$ is algebraically closed
in $\mathbb{C}(x,y,z)$, we obtain that $\mathbb{C}(x,y,z)^{\Phi}=\mathbb{C}(x)$. As a consequence, the ring of invariants of $\Phi$ is equal to $\CC[x]$.

\end{document}